\documentclass[11 pt, psamsfonts]{amsart}
\usepackage{amssymb}

\def\ignore#1{\relax}

\def\g{\mathfrak g}

\def\R{{\mathbb R}}
\def\Z{{\mathbb Z}}

\def\F{{\bf F}}
\def\la{\lambda}
\def\La{\Lambda}
 
\def\T{{\bf T}} 
\def\Pb{{\bf P}}

\def\W{\mathcal W}
\def\S{\mathcal S}
\def\H{\mathcal H}
\def\N{\mathcal N}
\def\A{\mathcal A}
\def\Ca{\mathcal C}

\def\one{\mathbf 1}
\def\Lm{{\bf F}_m}

\def\Ll{{\bf F}_\La}
\def\H{\mathcal H}

\def\Se{\mathcal S}
\def\A{\mathcal A}
\def\E{\mathcal E}
\def\T{{\bf T}}
\def\a{{\bf a}}
\def\b{{\bf b}}
\def\W{\mathcal W}
\def\M{\mathcal M}

\def\Nu{\underline{K}}

\def\ni{\noindent}

\def\v{\vskip 2.5mm}
\def\Lbar{\Z[v,v^{-1}]_0}
\def\inv{^{-1}}
\def\ignore#1{\relax}
\def\one{{\bf 1}}
\def\Nbar{\underline{N}}

{\theoremstyle{plain}
\newtheorem{theorem}{Theorem}[section]
}

{\theoremstyle{plain}

}

{\theoremstyle{plain}
\newtheorem{corollary}[theorem]{Corollary}
}

{\theoremstyle{plain}
\newtheorem{lemma}[theorem]{Lemma}
}

{\theoremstyle{plain}

}

{\theoremstyle{definition}
\newtheorem{definition}[theorem]{Definition}
}

{\theoremstyle{definition}

}

{\theoremstyle{exercise}

}

{\theoremstyle{remark}

}

{\theoremstyle{remark}

}

\numberwithin{equation}{section}

\renewcommand{\labelenumi}{{ \theenumi.}}
\renewcommand{\labelenumii}{{(\alph{enumii})}}

\ignore{
\input BoxedEPS
\SetRokickiEPSFSpecial
\HideDisplacementBoxes
\SetEPSFDirectory{./graphics/} }

\ignore{
\input BoxedEPS
\SetTexturesEPSFSpecial
\HideDisplacementBoxes
\SetEPSFDirectory{:graphics:}
}

\begin{document}

\title[A path algorithm for Kazhdan-Lusztig polynomials]
{A Path Algorithm for Affine Kazhdan-Lusztig Polynomials }

\author{Frederick M. Goodman and Hans Wenzl}
\thanks{H.W. was  supported in part by NSF grant \#DMS 9706839.}

\address{Department of Mathemematics\\ University of Iowa\\ Iowa City, Iowa}

\email{goodman@math.uiowa.edu}

\address{Department of Mathematics\\ University of California\\ San Diego,
California}

\email{wenzl@brauer.ucsd.edu}

\bigskip

\maketitle

Let $\la , \mu$ be dominant integral weights of a semisimple Lie algebra
$\g$ and let $l$ be a positive integer larger than the dual Coxeter number
of $\g$. Using these data and certain parabolic Kazhdan-Lusztig 
polynomials for affine Weyl groups, one can define polynomials 
$n_{\la ,\mu}$ in $\Z [v]$. It
has been shown recently that theses polynomials can be used to
compute the character of indecomposable tilting modules of quantum groups
at a root of unity (see [S1,2]).
For Lie type $A$, these polynomials also describe
the structure of the finite Hecke algebras of type $A$ for
$q$ a primitive $l$-th root of unity.
This was derived by a different method (see [A], ~\cite{Geck}, [LLT]),
motivated by Kashiwara's crystal bases;
the coincidence of these polynomials was shown in [GW2] and
~\cite{Varagnolo-Vasserot}.
The main result of our paper is an algorithm for computing the
polynomials $n_{\la ,\mu}$ using piecewise linear paths, for all Lie
types. This was, in part, inspired by the  algorithm of Lascoux,
Leclerc, and Thibon [LLT], which was used in the Hecke algebra setting. So
our algorithm can be viewed as a generalization of the LLT
algorithm to arbitrary weight lattices.
However, the method of proof here relies heavily on Kazhdan-Lusztig
theory.

One of the motivations for this work was to study in more detail
the behaviour of the polynomials  $n_{\la\mu}$ 
with $\la, \mu$ near the boundary of the Weyl chamber. As already
explained in our previous paper [GW2] for type $A$, the
algorithm is considerably faster than previous algorithms;
this allows a more extensive empirical investigation of the
polynomials  $n_{\la\mu}$.
Moreover, our algorithm shows some similarities with Littelmann's
path algorithm which describes the decomposition of tensor products
of Kac-Moody algebras. Indeed, at least in some special cases, our algorithm
can be used to compute the decomposition of the tensor product of
tilting modules into a direct sum of indecomposable tilting modules.
This could be used to compute the dimensions of simple modules of
the finite Hecke algebras of type $A$, using the $q$-analogue of
Schur-Weyl duality.
Moreover, for Lie type $A$, the path operators can be  identified 
with operators appearing in the Fock space
representation of $U_q\widehat{sl_l}$. It would be interesting if one
could find a similar interpretation of these operators also for other 
Lie types. This, as well as the 
connection to tilting modules and Hecke algebras
are discussed at the end of section 3.

Our paper is organized as follows. In section 1 we first review
basic facts of affine Hecke algebras and their Kazhdan-Lusztig polynomials.
We then study the multiplication of parabolic Kazhdan-Lusztig elements 
by Kazhdan-Lusztig elements $C_P$ corresponding to the longest element
of an arbitrary parabolic subgroup $P$. In Section 2, we define a 
geometric procedure which roughly corresponds to such multiplications.
In Section 3 we describe the already mentioned algorithm and discuss
various other (potential) applications.

$Acknowledgement$. Part of the work was done while both authors participated
at the workshop on representation theory of algebraic groups and quantum
groups in Aarhus. It is a pleasure to thank the organizers, Henning Hahr
Andersen and Jens Carsten Jantzen for hospitality and support.
\medskip
\section{Affine Hecke algebras}

We review basic facts of Kazhdan Lusztig polynomials for affine Hecke
algebras,  following Soergel [S1].

\subsection{}
Let $\Delta$ be a root system of a simple Lie algebra,
$Y$ its root lattice  and $X$ its weight lattice. Moreover,
let $V=X\otimes_\Z\R$. The highest root
is denoted by $\theta$. We fix a bilinear
form $\ \cdot\ $ on the $\Z$ span of $\Delta$ 
uniquely determined  by $\alpha\cdot\alpha=2$
for all short roots $\alpha$ (see [Lu, 1.1]).
We obtain an embedding of the root lattice into the weight lattice
via this bilinear form.

Define for each positive  root $\alpha$ the hyperplane
$E_{\alpha}\ =\ \{ x\in V,\ x\cdot\alpha=0\}$.
More generally, we define for any integer $r$ and any positive
root $\alpha$ the affine hyperplane
$$E_{\alpha ,r}= 
\{ x\in V,\ x\cdot\alpha = r\} .$$
Let $\W$ be the affine reflection group
generated by the finite Weyl group $W$ corresponding to $\Delta$ and
the affine reflection in the hyperplane $E_0= E_{\theta , l}$
where $\theta$  is the highest root  in $\Delta$. The image of the
defining action of an element $w\in \W$ on $x\in V$ is denoted by $w(x)$.
Let $\rho$ be half the sum of all positive roots in $\Delta$. Then
we define the dot action of $\W$ by the formula $w.x=w(x+\rho)-\rho$.

It can be shown that the collection $\E = \{ E_{\alpha ,ml}, 
\alpha\in \Delta_+, m\in\Z\}$ partitions $V$ into a disjoint collection
of fundamental domains, called alcoves, with respect to $\W$.
The hyperplane $E_\alpha=E_{\alpha ,0}$ divides $V$ into two (open)
halfspaces
$E_\alpha^+=\{ x\in V,\ (x,\alpha)>0\}$ and $E_{\alpha}^-$, the other
half space. The half spaces $E_{\alpha,m}^+$ and
$E_{\alpha,m}^-$ are the images of $E_\alpha^+$ resp $E_\alpha^-$ 
under a translation which maps $E_\alpha$ to $E_{\alpha ,m}$.
The dominant Weyl chamber $\Ca$ is equal to the intersection of 
$E_\alpha^+$, with $\alpha$ a simple root.
We say that a hyperplane $E\in\E$ contains a lower wall of the alcove
$A$ if $A$ is on the positive side of $E$, i.e. $A\subset E^+$
and the intersection of $E$ with the closure of $A$ 
has codimension 1;
similarly, we say that $E$ contains an upper wall of $A$ if
$A\subset E^-$ and the intersection of $E$ with the closure of $A$ 
has codimension 1.

The fundamental alcove $A_+$ is defined to be the intersection
of $\Ca$ with $E_{\theta ,l}^-$, where $\theta$ is the highest root.
The fundamental box $\Pi$ is the  set of points $x\in\Ca$ which 
satisfy $0<x\cdot \alpha\leq l$ for all
simple roots $\alpha$. The set of alcoves in the dominant Weyl chamber
is denoted by $\A_+$.
%

A facet $F$ is a nonempty connected subset of $V$ such that for any root 
$\alpha$ either $F$ is between two neighboring hyperplanes in $\E$
which are parallel to $E_\alpha$ (with empty intersection with
both of those hyperplanes), or it is contained
in $E_{\alpha ,ml}$ for some $m\in \Z$.
For a facet $F$ we define $A_+(F)$ to be the alcove which contains
$F$ in its boundary and which is on the positive side of each
hyperplane which contains $F$. The alcove $A_-(F)$ is defined similarly.

\subsection{Affine Hecke algebras}
The affine Hecke algebra $\H = \H(\W, \Se)$ is the associative algebra  with
identity element $\one$ over  the ring of Laurent polynomials
$\Z[v, v\inv]$ with generators $\{T_s : s \in \Se\}$ satisfying the  braid
relations and the quadratic relation $T_s^2 = v^{-2} \one + (v^{-2} -1) T_s$,
for $s \in \Se$.  Using the generators $H_s = v T_s$, one has instead
the quadratic relations $H_s^2 = \one + (v^{-1} - v) H_s$.  The Hecke algebra
has a basis $\{H_x : x \in \W\}$ satisfying $H_xH_y = H_{xy}$ in case
$\ell(xy) = \ell(x) + \ell(y)$, and $H_sH_x = H_{sx} + (v\inv-v)H_{x}$ in case
$\ell(sx) = \ell(x) -1$, for $x \in \W$ and $s \in \Se$; here $\ell$
denotes the usual length function on a Coxeter group with respect to
the generating set $\Se$.   The Hecke algebra has an involution
$d : a
\mapsto
\bar a$ defined by
$\bar v = v\inv$ and
$(H_x)^- = (H_{x\inv})\inv$.  An element fixed by this involution is called
self-dual.  We denote the self-dual part of $\Z[v,v\inv]$ by 
$\Z[v,v\inv]_0$.  Fundamental
self-dual elements are the
$C_s = H_s + v$ for
$s
\in \Se$.   Let $P$ be the parabolic subgroup generated by a subset
$\Se_P$ of
$\Se$.  We will 
assume here that $\S_P$ is a proper subset of $\Se$,  which implies
that $P$ is finite. We shall denote by $C_P$ the Kazhdan-Lusztig element
in the Hecke algebra $\H(P,\Se_P)$ corresponding to
the longest element $w_0$ of $P$. It is well-known that in our notations
$C_P$ is given by the formula
\begin{equation} \label{KL-element}
C_P=\sum_{w\in P} v^{\ell_P-\ell(w)}H_w,
\end{equation}
where $\ell_P=\ell(w_0)$. It is also well-known and easy to derive
that for any $w\in P$ we have $H_wC_P=v^{-\ell(w)}C_P$.
The $v$-order $[P]_v$ of $P$ is defined by
$$[P]_v=v^{-\ell_P}\ \sum_{w\in P} v^{2\ell(w)}= v^{\ell_P}
\ \sum_{w\in P} v^{-2\ell(w)}.$$
If $Q\subset P$ is a parabolic subgroup of $P$, and the element $C_Q$
is defined accordingly, we can easily derive from the formulas mentioned
before that 
\begin{equation} \label{multiplication of C_P's}
C_PC_Q=[Q]_vC_P=C_QC_P.
\end{equation}

\subsection{ Parabolic module} \label{parabolic}
Let $\Se_0 \subseteq \Se$ be the set of simple reflections fixing the origin. The
finite Weyl group $W$ is the Coxeter subgroup of $\W$ with generating set 
$\S_0$. Let $\H_f = \H(W, \Se_0) \subseteq \H$ denote its Hecke algebra.
Consider the sign representation of $\H_f \rightarrow \Z[v, v\inv]$ which
takes each $H_s$ to $-v\inv$.  Let $\N$ denote the induced right
$\H$-module
$$
\N = \Z[v, v\inv] \otimes_{\H_f} \H.
$$  Let $\W^f \subseteq \W$ be the coset  representatives of minimal length
of the right cosets of $W$ in $\W$.  Then $\N$ has a basis
$N_x = \one \otimes H_x$, for $x \in \W^f$, and the operation of  the $C_s$
for $s \in \Se$ on this basis has the following form:

\begin{equation} N_x C_s =
\begin{cases} N_{xs} + v N_x &\text{if $xs \in \W^f$ and $xs > x$}\\
N_{xs} + v\inv N_x &\text{if $xs \in \W^f$ and $xs < x$}\\ 0   &\text{if
$xs
\not\in \W^f$,}
\end{cases}
\end{equation}
 where the inequality signs refer to the Bruhat order on
$\W$.  The involution on $\H$ induces an involution on $\N$ defined by $a
\otimes b \mapsto \bar a \otimes \bar b$.

The following result is  discussed in 
~\cite{Soergel1}, following
~\cite{Deodhar1},~\cite{Deodhar2}.
\medskip

\begin{theorem} \label{Theorem KL basis} $\N$ has a unique basis
$\{\underline N_x : x
\in
\W^f\}$ satisfying
\begin{enumerate}
\item $\Nbar_x \in N_x + \sum_{y < x} v\Z[v] N_y$.
\item $\Nbar$ is self-dual.
\end{enumerate} 
\end{theorem}

\medskip
\v\ni
\begin{definition} The affine Kazhdan-Lusztig polynomials $n_{y,x}$ are
defined by
$$
\Nbar_x = N_x + \sum_{y < x} n_{y, x} N_y
$$
\end{definition}

\medskip The affine Weyl group $\W$ acts freely and transitively on alcoves,
so there is a bijection $\W \rightarrow \A$ given by $w \mapsto w A_+$, where
$A_+$ is the unique alcove in $\A_+$ containing 0 in its closure. Under this 
bijection, the elements of $\W^f$ correspond to alcoves contained in the 
positive Weyl chamber.  One also has an action of $\W$ on the right on $\A$
given by
$(w A_+) x = wx A_+$.  

Using the bijection between $\W^f$ and $\A_+$, one can  rename the
distinguished elements of the right $\H$ module
$\N$ using alcoves $A \in \A_+$ rather than coset representatives $x \in
\W^f$.   Thus if $x, y \in \W^f$ correspond to $A, B \in \A_+$, then we write
$N_A$ for $N_x$, $\Nbar_A$ for $\Nbar_x$, and $n_{A, B}$ for $n_{x, y}$. The
right action of $\H$ is then given by

\begin{equation} \label{equation right action of Cs} N_A C_s =
\begin{cases} N_{As} + v N_A &\text{if $As \in \A_+$ and $As \succ A$}\\
N_{As} + v\inv N_x &\text{if $As \in \A_+$ and $As \prec A$}\\ 0  
&\text{if $As \not\in \A^+$,}
\end{cases}
\end{equation} 
where now the inequalities have a geometric interpretation:
$As \succ A$ if $As$ is on the positive side of the hyperplane separating
the two alcoves. 
The elements $\Nbar_A$ are computed by a recursive
scheme as follows: One
has $\Nbar_{A^+} = N_{A_+}$.  Given $A \ne A_+$, one can choose $s \in \Se$
such that $As \in \A_+$ and $As \prec A$.  The element $\Nbar_A$ is then
computed by
$$\Nbar_A= \Nbar_{As}C_s - \sum_{B\prec A} f_{B,A}(0)\Nbar_B,$$
where $f_{B,A}$ is the coefficient of $B$ in $\Nbar_{As}C_s$.
\ignore{
$\Nbar_A$ one takes
$$
\Nbar_{As} C_s = N_A + \sum_{B \prec A} f_{B, A}(v) N_B.
$$ 
This element is self-dual, but may have coefficients with non-zero
constant term.  So one corrects these coefficients by subtracting a self-dual
linear combination of $\Nbar_B$ for $B \prec A$.}

\subsection{} \label{order}
Let $F$ be a facet. The left-stabilizer $\tilde P$ of $F$
is generated by the reflections in the hyperplanes containing $F$.
Let ${\mathcal A}_F$ be the set of all alcoves $A$ such that 
$F$ is in the boundary of $A$. Then it is well-known (see e.g. [Ja, 6.11])
that there exists a  1-1 correspondence
between the elements of $\tilde P$ and the elements of the  set 
${\mathcal A}_F$;
it can be defined by the map $w\in \tilde P\mapsto w(A_+(F))$. In particular,
$\tilde P$ is a finite reflection group, 
which is generated by a subset of reflections
in the hyperplanes containing a lower wall of $A_+(F)$.

Let $w_F\in \W$ be the element such that  $A_+(F)=w_F(A_+)$. Then
we define the right-stabilizer of $F$ to be the subgroup 
$P=w_F^{-1}\tilde P w_F$; one checks easily that it does
indeed coincide with all the elements
$w\in \W$ such that $\A_Fw=\A_F$, and hence $Fw=F$. It is generated
by a subset $\Se_F=\Se_P\subset \Se$ of simple reflections.
Then it follows from the definition of right action and the discussion above
that the map 
\begin{equation} \label{correspondence}
w\in P\mapsto A_-(F)w
\end{equation}
defines a 1-1 correspondence
between the elements of $P$ and the elements of the  set ${\mathcal A}_F$.
Moreover, if all the alcoves of $\A_F$ are contained in $\Ca$,
this correspondence is order preserving between the Bruhat order of $P$
and the order  $\prec$ on $\A_F$.
\vskip .2cm
Let $F_0$ be a facet which contains $F$ in its
closure. We define the quantities $\a_F(F_0)$ and $\b_F(F_0)$
to be the number of hyperplanes containing $F$ which are `above'
respectively `below' $F_0$; a hyperplane $H$ is said to be above $F_0$ if
$F_0$ is on the negative side of $H$, with a similar definition for
`below'.
We get the following formulas
which immediately follow from the definitions:
\begin{equation} \label{2.5.4}
\b_F(F_0)=\b_F(A_-(F_0))\quad {\rm and}\quad \a_F(F_0)=\a_F(A_+(F_0)).
\end{equation}
Also observe that if 
 $w\in P$ and $A=A_-(F)w$, then 

\begin{equation} \label{2.5.5} \b_F(A)=\ell(w)\quad {\rm and} \quad
\a_F(A)=\ell_P-\ell(w).
\end{equation}
\smallskip

\subsection{} \label{right action of NP}
We now consider the action of  the Hecke algebra $\H(\W, \Se)$ on 
the module $\N$. Let $\H(P,\Se_P)$ be the Hecke algebra
corresponding to the  right stablizer  $P$ of the facet $F$.
If all the alcoves of $\A_F$ are in $\Ca$, it follows from 
equation \ref{correspondence} and the text below it that the map

\begin{equation} \label{2.5.6}
H_w\in \H(P,\Se_P)\mapsto N_{A_-(F)w}\in \A_F
\end{equation}
defines an isomorphism of $\H(P,\Se_P)$ right modules.
The definition of the Kazhdan-Lusztig element $C_P$, see \ref{KL-element},
also suggests defining the element
\begin{equation} \label{2.5.7} 
N_F=\sum_{A\in \A_F} v^{\a_F(A)}N_A=N_{A_-(F)}C_P,
\end{equation}
where the last equality follows from \ref{2.5.5} and \ref{2.5.6}.
For a given facet $F$, 
we define a partial order $\prec$ on the set $\{N_{w(F)}, w\in \W$
such that $w(F)\subset \Ca\}$ by
\begin{equation} \label{facetorder}
F_1\prec F_2 \quad \leftrightarrow\quad  A_+(F_1)\prec A_+(F_2).
\end{equation}
Finally, let $F_0$ be a facet which contains $F$ in its boundary,
and let $Q\subset P$ be its stabilizer. Recall that $[Q]_v=v^{-\ell_Q}
\sum_{w\in Q}v^{2\ell(w)}$. It follows from the definitions of
$N_F$ and $N_{F_0}$ that
\begin{equation} \label{2.5.8}
N_F=\sum_{F_0'} v^{\a_F(F_0')}N_{F_0'},
\end{equation}
where the summation goes over the facets of the form $F_0w$ with $w\in P$.
The following elementary lemma is probably well-known to experts in 
Kazhdan-Lusztig theory.
\medskip
\begin{lemma} \label{rightmult}
Let $F_0$, $F$, $Q$ and $P$ be as just defined, with
$F$ in the closure of the dominant Weyl chamber $\Ca$, and let
$A\in \A_F$. Then we have
\begin{equation} \label{N_AC_P formula}
N_AC_P=\begin{cases} 0 & \text{if $F\not\subset \Ca$,}\\
v^{-\b_F(A)} N_F& \text{otherwise,}
\end{cases}
\end{equation}
and
\begin{equation} \label{N_{F_0}C_P formula}
N_{F_0}C_P=
\begin{cases} 0 & \text{if $F\not\subset \Ca$,}\\
{[Q]_v} v^{-\b_F(F_0)} N_F& \text{otherwise,}
\end{cases}
\end{equation}
\end{lemma}

\begin{proof}
Let $F$ be a facet in the boundary of $\Ca$ such that 
$F\not\subset \Ca$, and let $A$ be an alcove in $\A_F$.
We claim that $N_AC_P = 0$.
There exists an alcove  $B\in \A_F$, 
$B\subset\Ca$ which has a wall in a boundary hyperplane of $\Ca$.
Let $s$ be the simple reflection corresponding to this wall.
Then $N_BC_s=0$  (see e.g. [S1], or Section \ref{parabolic});
but then also $(v+v^{-1})N_BC_P= N_BC_sC_P=0$, which shows the claim
for $A=B$.   The claim will follow inductively if we can prove it for
any alcove $A \in \A_F$ which shares a wall with an alcove $B \in \A_F$
for which it has already been established.  Let  $s$ be the simple
reflection corresponding to such a  wall, and observe that $s \in P$.  
Then
$N_A=N_{Bs}= N_BC_s-v^{\pm 1}N_B$, from which one easily deduces the
claim for $A$.

Assume now that $F \subset \Ca$. Then necessarily any alcove $A\in \A_F$
has to be in $\Ca$.
If $A=A_-(F)$, we get $N_{A_-(F)}C_P=N_F$ from \ref{2.5.6} 
and the definitions of $N_F$ 
and of $C_P$, see \ref{2.5.7}  and
\ref{KL-element}. If $A=A_-(F)w$ for some $w\in P$,
then 
$$N_AC_P=N_{A_-(F)}H_wC_P=v^{-\ell(w)}N_{A_-(F)}C_P=v^{-\b_F(A)}N_F,$$
 using \ref{2.5.5} and \ref{2.5.6}.
To prove the second claim, observe that because of $Q\subset P$ we get
\ignore{
$$N_{F_0}C_P=\sum_{w\in Q} v^{\ell_Q-\ell(w)}N_{A_-(F_0)}H_wC_P=
\sum_{w\in Q} v^{\ell_Q-2\ell(w)} N_{A_-(F_0)}C_P= [Q]_vv^{-\b_F(F_0)}N_F,$$}
$$
N_{F_0}C_P = N_{A_-(F_0)} C_Q C_P = [Q]_v N_{A_-(F_0)}  C_P  = {[Q]_v} v^{-\b_F(F_0)} N_F,
$$
where the first equality follows from the definition of $N_{F_0}$ (formula \ref{2.5.7}), 
the second equality from formula \ref{multiplication of C_P's}, and the
last equality  from the formulas \ref{N_AC_P formula} for $N_A C_P$ and 
\ref{2.5.4}  for $\b_F(F_0)$.
\end{proof}

\subsection{}
\ignore{We keep the notations of previous subsections, with
$F_0$ being a facet with right stabilizer $Q$, and such that $F_0$ contains
the facet $F$ in its boundary.}
 We  define the self-dual element
$\Nbar_F$ by $\Nbar_F=\Nbar_{A_+(F)}$. Observe that for two facets $F$,
$\tilde F$, not necessarily of same dimension we have 
$\Nbar_F=\Nbar_{\tilde F}$ iff $A_+(F)=A_+(\tilde F)$.

\begin{lemma} \label{2.6}
\renewcommand{\labelenumi}{{(\alph{enumi})\quad}}
			\renewcommand{\labelenumii}{{(\roman{enumii})\ }}
Let $F$ be a facet with right stabilizer $P$, and let $F_0$ be a facet with right stabilizer $Q
\subseteq P$, such that $F_0$ contains the facet $F$ in its boundary.
Write $\Nbar_{F_0}=\sum_D r_D N_D$,
with $D\in \A_+$.  Then
\begin{enumerate}
\item
For any $A \in \A_F$, $\Nbar_AC_P$ is a  $\Lbar$ linear
combination of elements of the form $\Nbar_B$ with $B\prec A_+(F)$ (or
equal),
\item
 $\displaystyle \Nbar_{F_0}= \sum_{F_0'\in \W(F_0)}r_{F_0'} N_{F_0'}$, 
with $r_{F_0'}=r_{A_+(F_0')}$,

\item $\displaystyle \Nbar_{F_0}C_P=  [Q]_v\ \sum_{F_0' \in
\W(F_0)}r_{F_0'}v^{-\b_{F'}(F_0')}N_{F'},$ where $F'$ denotes, for each ${F_0'\in \W(F_0)}$,    
the unique facet conjugate to $F$ which is contained in the closure of $F_0'$. 
\item Let $F_1$ be any facet such that $F_1$ contains  $F$ in its boundary.
Then 
$$\Nbar_{F_0}C_P= [Q]_v\ \sum_{F_0'\in
\W(F_0)} \sum_{F_1'} r_{F_0'} v^{\a_{F'}(F_1')-\b_{F'}(F_0')}\ N_{F_1'},$$ 
where $F'$ is related
to $F_0$ as in part (c), and $F_1'$ ranges over all facets in $\W(F_1)$ which contain
$F'$ in their closure.
\end{enumerate}
\end{lemma}
 
\begin{proof}
To prove part (a), observe that $C_P$ is equal to the Kazhdan-Lusztig
element $C_{w_0}$ of the longest element $w_0$ in $P$. In particular,
it is self-dual and hence so is $\Nbar_AC_P$. Hence this product
can be written as a  $\Lbar$ linear combination of self-dual elements
$\Nbar_B$; in fact the same statement is true for $\Nbar_AC_w$
for any $w\in P$. We will show the second statement of the claim
more generally for $\Nbar_AC_w$ for any $w\in P$ by induction on the
length $\ell(w)$ of $w$.
If $\ell(w)=1$, $w=s$ is a simple reflection of $P$. In this case
the claim follows from the original proof by Kazhdan and Lusztig:
in fact, $\Nbar_AC_P$ is a linear combination of $\Nbar_B$ with
$B$ majorized by the larger of $A$ or $As$. Either of those elements
contains $F$ in its boundary, and hence is majorized by $A_+(F)$.
For the induction step, let $w=w's$, with $\ell(w')<\ell(w)$ and $s\in \Se_P$.
Then the claim follows by applying the induction assumption twice
for $\Nbar_AC_w=\Nbar_AC_{w'}C_s$.

Now let $A$ denote $A_+(F_0)$.  Let $\Se_Q$ be the set of simple reflections $s\in Q$. By definition
of
$A=A_+(F_0)$ and $Q$, we have $As \prec A$ for $s\in \Se_Q$. But then
$\Nbar_AC_s=(v+v\inv)\Nbar_A$ and
$\Nbar_AH_s=v\inv\Nbar_A$ for all $s\in S_Q$, using basic properties
of Kazhdan-Lusztig elements (see [KL]).
But then also $\Nbar_AH_w=v^{-\ell(w)}\Nbar_A$
for all $w\in Q$ from which one easily concludes
$\Nbar_AC_Q=[Q]_v\Nbar_A$.

Let $D$ be an alcove, and let $F_0'$ be the facet in the
boundary of $D$ which is conjugate to $F_0$. Then it follows from
Lemma \ref{rightmult} that $N_DC_Q$ is equal to a multiple of $N_{F_0'}$,
and hence also $\Nbar_{F_0}C_Q$ is a linear combination of $N_{F_0'}$'s.
On the other hand, we have $\Nbar_{F_0}C_Q=[Q]_v\Nbar_{F_0}$ by the
previous paragraph of this proof. Hence already $\Nbar_{F_0}$ itself
is a linear combination of  $N_{F_0'}$'s.
From this follows claim (b), as the coefficient of $N_{F_0'}$ is equal 
to $r_{A_+(F_0')}$.
Using (b) and Lemma \ref{rightmult},(b), we get
$$\Nbar_{F_0}C_P=\sum_{F_0'} r_{F_0'}N_{F_0'}C_P=
[Q]_v\ \sum_{F_0'}r_{F_0'}v^{-\b_{F'}(F_0')}N_{F'},$$
where $F'$ is the unique facet conjugate to $F$ which is contained in the
closure of $F_0'$. This proves part (c), and part (d)  follows from  $N_{F'}=
\sum_{F_1'} v^{\a_{F'}(F_1')}N_{F_1'}$ (see \ref{2.5.8}).
\end{proof}
\medskip

In our application of part (d) of this lemma, $F_1$ will be related to $F_0$ as follows:
$F_1$ will be
the facet ``on the other side of $F$ from $F_0$;" to be more precise, it is the facet reached if one
extends a line connecting a point in $F_0$ with a point in $F$ a little beyond $F$.
Observe that $F_1$  then lies in the affine subspace spanned by $F_0$ and has
the same dimension as $F_0$.

\medskip
\section{Path operators}
\subsection{} \label{KL for points}
We now define Kazhdan-Lusztig polynomials $n_{\la\mu}$
for {\it arbitrary points} $\la$, $\mu$ in the dominant Weyl chamber 
$\Ca$ of $\g$ as follows: 
Let $\mu$ be a point in $\Ca$, and
let $F$ be the facet
containing $\mu+\rho$. If $\mu+\rho$ does not lie on any facet,
we define both $F$ and $A_+(\mu)$ to be the alcove
containing  $\mu+\rho$; otherwise, we define $A_+(\mu)$ to be equal
to $A_+(F)$.
Then the polynomials  $n_{\la\mu}$ are defined by
\begin{equation}
n_{\la\mu}=
\begin{cases} n_{A_+(\la),A_+(\mu)} & \text{if $\la+\rho\in \W(\mu+\rho)$,}\\ 
0& \text{otherwise.}
\end{cases}
\end{equation}

The elements $\Nbar_\mu$ are defined in the free $\Z[v,v^{-1}]$ module
${\mathcal M}$ with basis $\{ N_\la , \la\in \Ca\}$ by
\begin{equation}\label{Nbar}
\Nbar_\mu = \sum_\la n_{\la\mu}N_\la .
\end{equation}
By definition, the polynomials $n_{w.\mu ,\mu}$ do not change if we vary
$\mu$ within its facet, and keep $w\in\W$ fixed.
If we denote by $F(\la+\rho)$ the facet containing $\la+\rho$,
and if $\mu+\rho$ is on the facet $F_0$, we get
\begin{equation} \label{2.7.1}
\Nbar_{F_0}=\Nbar_{A_+(F_0)}=\sum_{\la} n_{\la\mu}N_{F(\la+\rho)}.
\end{equation}
\ignore{
It is shown in [S1,2] that the polynomials $n_{\la\mu}$ describe the
Weyl filtration of the indecomposable tilting module with highest weight
$\mu$, for $\mu$ a dominant integral weight.
}
An element $M$ in ${\mathcal M}$ is said to be self-dual 
if it can be rewritten as a
linear combination of elements of the form $\Nbar_\mu$ with coefficients
in $\Lbar$. If $M = \sum_{\la} m_{\la}N_\la$, it follows
that $M$ is self-dual if and only if the element
$\sum_{\la} m_{\la}N_{F(\la+\rho)}$
is self-dual.
\v\ni

\ignore{
{\bf 2.6 Lemma} Assume that $\Nu_\mu$ and $\Nbar_\mu$ coincide. Moreover,
let $\La=t\La_m$, where $\La_m$ is a fundamental weigth, and $t$ is a
positive scalar $\leq 1$ such that the line segment connecting 
$\mu+\rho$ with $\mu+\rho+\La$ intersects at most one facet $F$
with smaller dimension. Then 
$\Lm\Nbar_m$ is a multiple of $\Nbar C_P$
where $P$ is the stabilizer of $F$.
\v
$Proof.$ Follows from the definitions.
}

\subsection{} \label{geometric}
We describe a geometric procedure which will be useful
for describing the behaviour of Kazhdan-Lusztig elements $\Nbar_\mu$ under
multiplication by $C_P$.
Let $F_0$ be the facet containing $\mu+\rho$, and let $F$ be a facet in
the boundary of $F_0$. Let $\Lambda$ be a vector such that the line
segment between $\mu+\rho$ and 
$\mu+\rho+\Lambda$ is contained in $F_0\cup F\cup F_1$, with $F_1$
the facet containing $\mu+\rho+\Lambda$. In applications we will
also encounter the 
degenerate cases with $\mu+\rho\in F$ (i.e. we have a path going
from a facet $F$ of smaller dimension to a facet $F_1$ which contains
$F$ in its boundary), and the case with $\mu+\rho+\La\in F$
(i.e. we have a path going from the facet $F_0$ to a facet in its 
boundary). Obviously, the general case described first can be obtained
as a combination of the two degenerate cases. In the following
we will only deal with the general case; the degenerate cases can
be treated similarly.

Let $\a_F(\la)=\a_F(F_0)$,
where $F_0$ is the facet containing $\la+\rho$ (see \ref{2.5.4}),
with the definition of $\b_F(\la)$ similar. Then
we define
\begin{equation}\label{def1} 
\Ll N_\mu=
\begin{cases} \ 0
& \text{if $F\not\subset \Ca$,}\\ 
\displaystyle\ \sum_{\gamma\,  :\,  {\gamma +\rho\in (\la+\rho+\Lambda)P}} 
v^{\a_F(\gamma)-\b_F(\la)}N_\gamma
& \text{if $F\subset \Ca$.}
\end{cases}
\end{equation}
We extend the definition of $\Ll$ to the elements $\Nbar_\mu$ in
the following nontrivial way: Observe that $\Nbar_\mu$ is a  linear
combination of elements $N_\la$ with $\la\in \W.\mu$. By assumption on
$\Lambda$, the line segment  between $\mu+\rho$ and 
$\mu+\rho+\Lambda$ is in the affine subspace spanned by $F_0$.
Hence if $w.\mu=\la$, we obtain a well-defined line segment
from $w(\mu+\rho)$ to $w(\mu+\rho+\Lambda)$, independent of the
choice of $w$. We denote the vector given
by it by $\bar w(\Lambda)$, where $\bar w$ is the image of $w$ in the
quotient $\W/T\cong W$, with $T$ being the normal subgroup of translations
in $\W$. If $\Nbar_\mu = \sum_{\la\in \W.\mu} n_{\la\mu}
N_\la$, we define
\begin{equation}\label{def2}
\Ll\Nbar_\mu = \sum_{\la=w.\mu} n_{\la\mu}{\bf F}_{\bar w(\Lambda)}
N_\la.
\end{equation}

\begin{lemma} \label{2.8}
The element $\Ll \Nbar_\mu=\sum_\gamma u_\gamma N_\gamma $ 
is self-dual. Moreover, if $Q$ is the stabilizer of the facet $F(\mu+\rho)$
containing $\mu+\rho$, then $u_\gamma$ equals the coefficient of
$N_{F(\gamma +\rho)}$ in the expansion of $[Q]_v\inv\Nbar_{F(\mu+\rho)}C_P$ as in
Lemma \ref{2.6}(d).
\end{lemma}
\v
\begin{proof}
Let $\la =w.\mu$, and let $F_0', F'$ and $F_1'$ be the facets
obtained by applying $w$ to $F_0, F$ and $F_1$. Observe that $\la+\rho\in
F_0'$ and that the line segment from $\la+\rho$ to $\la+\rho +\bar w(\Lambda)$
is contained in $F_0'\cup F'\cup F_1'$. Moreover
$F_0'$ and $F'$ have right stabilizers $Q$ and $P$ respectively.
It follows from Lemma \ref{2.6}(c) that
$$N_{F_0'}C_P=[Q]_v
v^{-\b_{F'}(F_0')}N_{F'}=[Q]_v\sum_{F_1'} v^{\a_{F'}(F_1')-\b_{F'}(F_0')}
N_{F_1'},$$
where the summation goes over all facets $F_1'$ conjugate to $F_1$ which
contain $F'$ in their boundary. 
Using the  1-1 correspondence between the facets $F_1'$ in the sum above,
and the points $\gamma + \rho$
in the orbit of $\la +\rho +\Lambda$ under the right action of $P$,
we obtain similarly
$${\bf F}_{\bar w(\Lambda)}N_\la =
\sum_{\gamma} v^{\a_{F'}{(\gamma)}-\b_{F'}(\la)}
N_\gamma
=\sum_{\gamma +\rho\in F_1'} v^{\a_{F'}(F_1')-\b_{F'}(F_0')}
N_{\gamma}.$$
The claim about the coefficient of $N_{F(\gamma +\rho)}=N_\gamma$ follows
from the last two equations, and the definition of $\Ll$, see \ref{def1}
and \ref{def2}. As $C_P$ is self-dual, so is $\Nbar_{F(\mu+\rho)}C_P$.
From this follows that $\Ll\Nbar_\mu$ is self-dual, by the last
two equations and the definition of duality on the module $\M$ (see end of
section
\ref{KL for points}).
\end{proof}

\subsection{} \label{positive}
In order to compute the $\Nbar_\mu$ explicitly, 
we define  a `positivity' operator $\Pb$  on self-dual linear combinations
$R=\sum_{\la \leq \mu} a_\la N_\la$, with $a_\la\in \Z[v,v^{-1}]$
and $a_\mu = 1$ as follows:
Let  $i_o$ be the smallest exponent of $v$ occurring in any of the
coefficients $a_\la$ with $\la <\mu$. If $i_0> 0$ we define $\Pb R=R$.
If $i_0\le 0$, we get rid of the lowest exponent by
the operation
$$R \quad \mapsto\quad R\ -\ \sum_{\la <\mu} 
(v^{i_o}+v^{-i_o})[(v^{-i_o}a_\la)(0)]\Nbar_\la;$$
if $i_0=0$ we subtract from $R$ the expression 
$\sum_{\la <\mu} a_\la(0)\Nbar_\la$.
Iterating this operation, as long as the lowest exponent is nonpositive,
we finally obtain a linear combination $\Pb R =\sum_\la b_\la N_\la$ with
$b_\la\in v\Z[v]$ for all $\la<\mu$, i.e. we obtain  (cf. Theorem
\ref{Theorem KL basis})
\begin{equation} \label{2.9.1}
\Pb (R) =   \Nbar_\mu .
\end{equation}

\begin{theorem} \label{2.9}
$\Nbar_{\mu+\Lambda}=\Pb\Ll\Nbar_\mu$ if $\Lambda$ is
in $\Ca$ satisfying the conditions at the beginning of Section 2.2.
\end{theorem}

\begin{proof}
 It was shown in Lemma \ref{2.8} that $\Ll \Nbar_\mu$ is a $\Lbar$-
linear combination of $\Nbar_\gamma$s. Moreover, using  the fact that
$\Lambda$ is in $\Ca$, we check easily that $\Nbar_{\mu +\Lambda}$ has
coefficient 1 in $\Ll \Nbar_\mu$. 
The claim follows from this and \ref{2.9.1}.
\end{proof}

\medskip
\section{ A fast algorithm}

\subsection{}
We now outline how our path version of the
 Kazhdan-Lusztig algorithm can be used to compute the coefficients
$n_{\la\mu}$ in an efficient way. Here,
the strategy is to stay at points with stabilizer
as large as possible as long as possible.

Let $\La$ be a vector in $\Ca$. Then we define the operator
$\T_\La :\M\to\M$ by $\T_\La \Nbar_\mu = \Nbar_{\mu+\La}$. 
If for given $\mu$ the vector $\La$ is small enough that 
the line segment between $\mu+\rho$ and  $\mu+\La+\rho$ satisfies
the conditions in Section \ref{geometric},  
$\T_\La \Nbar_\mu=\Pb\Ll\Nbar_\mu$,
by Theorem \ref{2.9}. For the general case, we choose a piecewise linear
path with vertices $\mu_0 +\rho = \mu+\rho, \mu_1+\rho,..., 
\mu_N+\rho=\mu+\La+\rho$ such that each directed line segment from
$\mu_i+\rho$ to $\mu_{i+1}+\rho$  lies in $\Ca$ and satisfies the conditions at the
beginning of section \ref{geometric}. This can always be achieved if
$\La\in \Ca$; in this case we can write $\La$ as a linear combination
of the fundamental weights with non-negative scalars.
So it suffices to chose as line segments 
small enough multiples of the fundamental weights.
Theorem \ref{2.9} then implies

\begin{corollary}
 $\T_\La \Nbar_\mu = (\Pb\F_{\mu_N-\mu_{N-1}})
(\Pb\F_{\mu_{N-1}-\mu_{N-2}})\ ...\  (\Pb\F_{\mu_1-\mu_0})\Nbar_\mu$,
independent of the choice of the piecewise linear  path.
\end{corollary}

\subsection{} \label{algorithm}
A dominant weight $\mu$  is called  {\it critical} if 
$(\mu +\rho)\cdot \alpha$
is divisible by $l$ for all roots $\alpha$ of $\g$.
In the simply-laced case, the smallest critical weight is
the Steinberg weight $(l-1)\rho$; in the nonsimply-laced case
the smallest critical weight is a multiple of $\rho$ or
$\check \rho$ (half the sum of the co-roots) depending on whether
$l$ is divisible by the ratio of the square lengths of a long
and a short root. In the following, we shall assume $\g$ to
be {\it simply-laced} in order to avoid having to deal with
various cases.
We call a dominant integral weight $\mu$ {\it interior} if 
$\mu\cdot \alpha_i \geq l$ for all simple roots $\alpha_i$.
The {\it fundamental box} $\Pi$ is the set
$$
\{ \sum_{i = 1}^{k-1}  t_i \Lambda_i : 0 < t_i \le l \}.
$$ 
The dominant Weyl chamber
is tiled by translates of the fundamental box $\Pi$ of the
form $\Pi + \mu_c+\rho$,
where $\mu_c$ is a critical weight in $ \bar\Ca$.
In particular, for each interior weight $\mu$, there is a unique
critical weight $\mu_c\in \Ca$ such that
$\mu +\rho \in \Pi+ \mu_c + \rho$, or, equivalently,
with $\mu \in \mu_c + \Pi$.

{\it  Case 1: Interior weights} 
Assume that the critical point $\mu_c$ of $\mu+\rho$ is in 
the interior. Then it is known that $\Nbar_{\mu_c}=N_{\mu_c}$ (see
~\cite{Soergel1}, ~\cite{Deodhar1}). We 
can compute $\Nbar_\mu=\T_{\mu-c_\mu}N_{\mu_c}$, with $\T_{\mu-c_\mu}$
realized by a suitable piecewise  linear path (see subsection 3.1).

{\it Case 2: Noninterior weights}
If $\mu$ is a non-interior weight, then $\mu + \rho =\sum_i a_i\La_i$,
with the $a_i$ non-negative integers with at least one of them $< l$.
We define scalars $b_i$ and $c_i$ as follows:
If $a_i<l$, then $b_i=c_i=a_i$; otherwise we set $c_i=l [a_i/l]$
and $b_i=l$, where for the purpose of this definition $[x]$ means the
integer part of the real number $x$. We then define dominant weights
$\mu_0$ and $\mu_1$ by
$$\mu_0+\rho =\sum_i b_i\La_i\quad {\rm and}\quad 
\mu_1+\rho =\sum_i c_i\La_i;$$
observe that $\mu_0$ is in the fundamental box, while $\mu_1$
is close to $\mu$, playing the role of the critical weight $\mu_c$ in case 1.
The computation of $\Nbar_\mu$ now goes in three steps:

(i) Compute $\Nbar_{\mu_0}$ directly via a suitable piecewise linear
path from $\rho$ to $\mu_0+\rho$.

(ii) Compute $\Nbar_{\mu_1}$ via a suitable piecewise linear path
from $\mu_0+\rho$ to $\mu_1+\rho$; suitable here means that the
path only goes through facets whose dimensions are at most one
higher than the dimension of the facet containing $\mu_1+\rho$.

(iii) Compute  $\Nbar_{\mu}$ via a path  from  $\mu_1+\rho$ to
$\mu+\rho$.

\subsection{} The main savings of this algorithm
over the usual algorithms comes from
step (ii) in the case of noninterior weights, for two reasons:
Firstly, working with facets requires much less memory
than having to work with all the alcoves surrounding it separately,
and hence also less operations for computing a new Kazhdan-Lusztig
element. Secondly, going through a facet of large codimension
requires less wall-crossings than having to go around via
simple wall-crossings.
Indeed, computer experiments for type $A$ show that this algorithm
produces enormous improvements in efficiency over the original LLT
algorithm, or the Soergel algorithm (see [GW2] for details).

\subsection{} Let $\la$ be a Young diagram with $< k$ rows. We identify
it with a weight of $sl_k$ as usual.  Taking for $\Lambda$ a 
weight of the generating
$k$-dimensional representation $V$ of $sl_k$, it is not hard to check
that the action of the operator $\Ll$ on $\Nbar_\la$
corresponds to the action of one of the
operators $f_i$ on $\la$, where the $f_i$'s
describe  the action of a Borel algebra of 
$U_q\widehat{sl_l}$ on the Fock space
consisting of all Young diagrams, see ~\cite{Misra-Miwa}. 
Indeed, this action served as motivation for
work by ~\cite{LLT},  ~\cite{Ariki} and also for this paper.
It would be interesting to find out whether one can find similar
interpretations for suitably chosen path operators, as defined in
this paper, for other weight lattices.

\subsection{Littelmann path algorithm} 
The
path operators in this paper can be used to describe the decomposition
of the tensor product of an indecomposable tilting module $T_\mu$
with highest weight $\mu$ (see ~\cite{Andersen tilting}) with the
fundamental representation $V$ of $U_qsl_k$ (as in the previous subsection).
Indeed, for given dominant integral weight $\mu$, 
one can find a subset $\{ \varepsilon_i\}_i$ of
the weights of $V$ such that the application of the operator 
$\sum_i {\bf F}_{\varepsilon_i}$ results into adding $all$ possible
weights of $V$ to $\mu+\rho$; e.g. if $\mu+\rho$ does not lie on any
hyperplane, one would have to take as subset $all$ weights of $V$.
Writing $(\sum_i {\bf F}_{\varepsilon_i})\Nbar_\mu=\sum_\nu m_\nu \Nbar_\nu$,
with $m_\nu\in \Lbar$ (which is always possible), it follows that the
multiplicity of the indecomposable tilting module $T_\nu$ with highest
weight $\nu$ in $T_\mu\otimes V$ is equal to $m_\nu(1)$. This is
an easy consequence of the Littelmann algorithm (see ~\cite{Li}) 
in the generic case,
and the fact that the multiplicity of an indecomposable tilting module
in the tilting module $T$ is completely determined by the character of $T$.
Hence, in this case, our algorithm
can be seen as an analog of the Littelmann path algorithm
for tilting modules (in a rather trivial special case).
The same reasoning will also work for other Lie types if one takes
for $V$ a minuscule representation.
It would be interesting
to see whether one could extend this observation to a description
via paths of tensoring $T_\la$ with an arbitrary tilting module $T$.
Note, however, that even if $T$ is a simple tilting module,
one can not adapt Littelmann's
algorithm in a straightforward way in general.

Finally we would like to recall that
one obtains the dimension of the simple module
$D^\mu$ of a Hecke algebra of  type $A$ 
(see ~\cite{Dipper-James1}) as the multiplicity of the
tilting module $T_\mu$ in $V^{\otimes n}$. 
This is a direct consequence of the $q$-analogue
of Schur duality between Hecke algebras of type $A$ and $U_qsl_k$
(see  ~\cite{Donkinb}, ~\cite{DPS}).
So the procedure
described in the previous paragraph is a method to compute such
dimensions inductively. Inductive procedures of this type have
already been obtained before by Kleshchev for symmetric groups
in positive characteristic and by Brundan (~\cite{Brundan}) for Hecke 
algebras by different methods.

\end{document}